\documentclass[review]{elsarticle}

\usepackage{amsmath, amssymb} % define this before the line numbering.
\newcommand\norm[1]{\left\lVert#1\right\rVert}
\usepackage{gensymb}
\usepackage{algorithm}
\usepackage{algpseudocode}
\usepackage{lscape}
\usepackage{subfigure}
\usepackage{amsthm}
\usepackage{amssymb}

\usepackage{tikz}
\usetikzlibrary{arrows}
\usepackage{verbatim}
\usetikzlibrary{positioning,fit,calc}
\usetikzlibrary{arrows}
\usepackage{tabularx}
\graphicspath{ {./images/} }
\DeclareGraphicsExtensions{.pdf,.png,.jpg}
\usepackage[font=small,labelfont=bf,tableposition=top]{caption}
\usepackage{lipsum}
\usepackage{multirow}
\usepackage{color}
\usepackage{microtype}
\usepackage{neuralnetwork}
\usepackage{xpatch}
\usepackage{mathrsfs}
\usepackage{mathtools}
\DeclarePairedDelimiter\ceil{\lceil}{\rceil}

\journal{Journal of \LaTeX\ Templates}

\begin{document}

\begin{frontmatter}

\title{On the approximation of a matrix}

\author{Samriddha Sanyal} \corref{mycorrespondingauthor}
\cortext[mycorrespondingauthor]{Corresponding author}
\address{Indian Statistical Institute, 203 B.T. road, Kolkata 700108, India}
\ead{samriddha.s@gmail.com}

% %% or include affiliations in footnotes:
% \author[mymainaddress,mysecondaryaddress]{Elsevier Inc}

% \author[mysecondaryaddress]{Global Customer Service\corref{mycorrespondingauthor}}
% \cortext[mycorrespondingauthor]{Corresponding author}

% \address[mymainaddress]{1600 John F Kennedy Boulevard, Philadelphia}
% \address[mysecondaryaddress]{360 Park Avenue South, New York}

\begin{abstract}
Let $F^{*}$ be an approximation of a given $(a \times b)$ matrix $F$ derived by methods which are not randomized. We prove that for a given $F$ and $F^{*}$,  $H$ and $T$ can be computed by randomized algorithm such that $(HT)$ is an approximation of $F$  better than $F^{*}$.

\end{abstract}


\end{frontmatter}

%\linenumbers

\section{Introduction}

Throughout the paper, we denote by $\norm{F}$ the Frobenius norm of a matrix $F$ i.e. the square root of the sum of the squares of the elements of $F$. For the desired rank $r$, the low-rank matrix approximation problem is to find $F_{apx}$ such that:
\begin{equation}
F_{apx}=\underset{rank(\hat{F}) \leq r}{\min} \norm{F- \hat{F}}.
    \label{prb}
\end{equation}
 The methods to approximate a given matrix can be classified into two groups that play a rather important role and have attracted a lot of interest: The first group is the group of all algorithms that are not randomized. The first group consists of sparsification \cite{achlioptas2007fast,gittens2009error}, column selection methods \cite{drineas2006fast, frieze2004fast}, approximation by dimensionality reduction \cite{papadimitriou2000latent,clarkson2009numerical} and approximation by submatrices  \cite{mahoney2009cur, goreinov1997theory}. The second group is the randomized algorithms \cite{halko2011finding, martinsson2016randomized} for low-rank matrix approximation.

Recently the randomized algorithms get more attention because of their better performance than other low-rank matrix approximation methods. However, so far, no direct theoretical superiority of the randomized algorithms over the other matrix approximation methods has been established.

In this paper, we prove the superiority of randomized algorithms over the other low-rank matrix approximation methods. We consider an approximation of $F$ say $F^{*}$ which is derived by any methods that are non-randomized with approximation error $\epsilon$. i.e. $\norm{F - F^{*}}= \epsilon$. We show that for a given approximation $F^{*}$ of $F$ with approximation error $\epsilon$, we can derive $H$ and $T$ by randomized algorithm such that $HT$ is an approximation of $F$ and $\norm{F - HT} < \epsilon$. That implies $HT$ is an approximation of $F$ better than $F^{*}$. Next, we discuss the result and expound the proof in detail.

\section{The result}
The goal of our paper is to show the efficacy of randomized algorithm for low-rank matrix approximation over the algorithms which are not randomized. This problem can be stated formally in the following way.

\newtheorem{thm}{Theorem}
\begin{thm}
Let $F^{*}$ be an arbitrary approximation of $F$ derived by the non-randomized algorithms for low-rank matrix approximation. The approximation error of $F^{*}$ is $\epsilon$. i.e. $\norm{F-F^{*}}= \epsilon$. Then for any such $F^{*}$, a $H$ and $T$ can be found by randomized algorithm such that $(HT)$ is an approximation of $F$ with approximation error $\norm{F-HT} < \epsilon$.
%\vspace{-0.15in}
\label{l1}
\end{thm}

\begin{proof}
%   Let we have the groundtruth $F$ that registers the players of the video frame into the top view model. Given a video frame, the problem is to approximate $F$. Let, the state-of-the-art method approximates $F$ as $F^{*}$ defined from (b) into (d). $F^{*}$ registers players from the video frame into the top-view model with approximation error $\epsilon$ i.e. $\norm{F-F^{*}} = \epsilon$, where $\norm{.}$ is the Frobenius matrix norm. Given arbitrary $F^{*}$, we aim to construct $T$ and $H$ such that $(HT)$ is a better approximation of $F$. More precisely, given arbitrary $F^{*}$, there exists $T$ and $H$ such that  $(HT)$ has approximation error less than $\epsilon$. i.e. $\norm{F- (HT)} < \epsilon$. 
  
  The proof is comprised of three steps. First, we construct $H$ and $T$ by randomized algorithms. The construction of $H$ and $T$ depends on the value of a variable $s$. Thereafter, we derive the error bound for the proposed approximation $(HT)$ involving $s$. Finally, the value of $s$ is computed for the given $\epsilon$ and using the value of $s$, $H$ and $T$ can be computed specifically with approximation error less than $\epsilon$. Next we expound the details.

%Assume $F$ is a matrix of dimension $(a \times b)$.
%$H$ is a homography that establishes correspondences between the global view of the ground and the top-view as shown in Fig. \ref{block}. Therefore, $H$ is invertible. So, essentially, we have to approximate $F$. 

\subsection{Construction of $H$ and $T$}

The task of computing a low-rank approximation to a given matrix $F$ can be divided into two steps. The first step is to construct a low dimensional subspace that captures the range of $F$. The second is to restrict the $F$ to the subspace and then compute a standard factorization (for example QR factorization) of the reduced $F$. Next, we discuss how the randomized algorithms handle the steps to compute low rank matrix approximation.

\begin{algorithm}[b]
   % \centering
    \caption{ Approximation by randomized algorithm (Proto-algorithm)}
Given $(a \times b)$ matrix $F$, target rank $r$ and an oversampling parameter $s$ this process computes an $a \times (r+s)$ matrix $H$ whose columns are orthonormal and whose range approximates the range of $F$.
\begin{enumerate}
 \item Generate $b \times (r+s)$ matrix $\mathscr{G}$.
  \item Compute $(F\mathscr{G})$.
  \item Construct a matrix $H$ whose columns form an orthonormal basis of the range $(F\mathscr{G})$.
\end{enumerate}
	\label{a1}
\end{algorithm}

%The first step can be executed with random sampling methods \cite{martinsson2016randomized}. 

To understand how randomness works, it is helpful to consider a motivating example \cite{martinsson2016randomized}. Let us consider $F= B + E$, where $B$ captures the range of $F$ and $E$ is a small perturbation error. Let the dimension of $F$, $B$ and $E$ be $(a\times b)$. Our aim is to obtain a basis of exact rank $r$ that covers as much of the range of $B$ as possible. In order to obtain $r$ rank approximation of $B$, a small number $s$ is fixed and $(r+s)$ random vectors $\{F(\alpha_{i})\}^{(r+s)}_{i=1}$ are generated:

\begin{equation}
     F(\alpha_{i})= B(\alpha_{i}) + E(\alpha_{i}),
     \label{factor1}
 \end{equation}
for $i=1,..., (r+s)$. The perturbation $E$ deviates the direction of each $F(\alpha_{i})$ outside the range of $B$. As a result the span of only $r$ vectors $\{F(\alpha_{i})\}^{r}_{i=1}$ may not cover the entire range of $F$. In contrast, the enriched set $\{F(\alpha_{i})\}^{(r+s)}_{i=1}$ enhance the chance of spanning the required subspace. The intuitive approach is applied in randomized algorithm for low rank matrix approximation. The randomized algorithm computes low rank approximation of $F$ in three steps as follows \cite{halko2011finding}:

First, a random $(b \times (r+s))$ matrix $\mathscr{G}$ is generated whose columns are Gaussian vectors. Thereafter compute $(F\mathscr{G})$. Finally, construct a matrix $H$ whose columns form an orthonormal basis of the range $(F\mathscr{G})$. The process is described in algorithm \ref{a1}. Once we get the $H$, then we can compute the other factor $(H^{*}F)$. Considering $T=H^{*}F$, the $F$
is approximated in factor form as $HT$.i.e. $F \approx HT$.

The approximation $HT$ is a variable of the oversampling parameter $s$. Now our objective is to compute the value of $s$ for which the $\norm{F-HT} < \epsilon$. Next, we compute the bound of expected approximation error of $\norm{F-HT}$. There after, we determine the value of $s$ that is needed to compute $H$ and $T$ so that the approximation error is less than $\epsilon$.
%$\Sigma_{1}= diag(\sigma_1,..., \sigma_r)$

\subsection{Computing the error bound}

\newtheorem{lemma}{Lemma}
\begin{lemma}
We aim to show: 
\begin{equation*}
    E(\norm{F-HT})) \leq (1+ \frac{r}{s-1}) ( \sum_{i=r+1}^{\min{(a,b)}} \sigma_{i}^2)
\end{equation*}
where $E$ is the expectation, $\sum_{i=r+1}^{\min{(a,b)}} \sigma_{i}^2$ is the theoretically minimal error in approximating F by a matrix of rank $r$.
%\vspace{-0.15in}
\label{l2}
\end{lemma}

\begin{proof}

First, consider the singular value decomposition of $F$ as $F= U_{1} \Sigma_{1} V_{1}^{*}$, where $U_{1}$ is a $(a \times r)$ orthonormal matrix, $\Sigma_{1}$ is a diagonal matrix containing the non negative singular values of $F$ and $V_{1}$ is a $(r \times n)$ orthonormal matrix. We call $U_{1}$ and $V_{1}$ as left unitary factor and right unitary factor respectively.
First partition the $\Sigma_{1}=[\Sigma_{2}|\Sigma_{3}]$, where the $\Sigma_{2}$ and $\Sigma_{3}$ are the diagonal matrix containing the first $r$ and $(b-r)$ singular values respectively. Thereafter, partition $V_{1} = [V_{2} | V_{3}]$ into blocks containing $r$ and $b - r$ columns respectively. Define $\mathscr{G}_{2} = V^{*}_{2}\mathscr{G}$ and $\mathscr{G}_{3} = V^{*}_{3}\mathscr{G}$. Since, $V_{2}$ and $V_{3}$ are orthonormal, then $\mathscr{G}_{2}$ and $\mathscr{G}_{3}$ are also Gaussian. We denote the pseudoinverse of $\mathscr{G}_{2}$ and $\mathscr{G}_{3}$ as $\hat{\mathscr{G}_{2}}$ and $\hat{\mathscr{G}_{3}}$
respectively. $\mathscr{G}_{2}$ and $\mathscr{G}_{3}$ are non overlapping, so they are stochastically independent. Applying Holder's inequality, we can write: 

\begin{equation}
E(\norm{F-H(H^{*}F)}))\leq (E(\norm{F-H(H^{*}F)}^{2}))^{1/2}
\label{f2}
\end{equation}

It is proved in \cite{halko2011finding} that:
\begin{equation}
   E(\norm{F-H(H^{*}F)}^{2}) \leq (\norm{\Sigma_{3}}^{2}_{F} + E(\norm{\Sigma_{3}\mathscr{G}_{3} \hat{\mathscr{G}_{2}})}^{2})
   \label{f3}
\end{equation}

 Therefore, using Eq. \ref{f2} and \ref{f3}, we can write:
 \begin{equation}
     E(\norm{F-H(H^{*}F)}))\leq (\norm{\Sigma_{3}}^{2}_{F} + E(\norm{\Sigma_{3}\mathscr{G}_{3} \hat{\mathscr{G}_{2}})}^{2})^{1/2}
     \label{f4}
 \end{equation}

We compute $E(\norm{\Sigma_{3}\mathscr{G}_{3} \hat{\mathscr{G}_{2}}}^{2})$ by conditioning on the value of $\mathscr{G}_{2}$ as follows:

\begin{equation}
     E(\norm{\Sigma_{3}\mathscr{G}_{3} \hat{\mathscr{G}_{2}}}^{2})   = E(E(\norm{\Sigma_{3}\mathscr{G}_{3} \hat{\mathscr{G}_{2}}}^{2})|\mathscr{G}_{2}) 
     \end{equation}

 The Frobenious norm is unitarily invariant. i.e. for any two orthonormal matrices $U_{1}$ and $V_{1}$, we can write $\norm{U_{1}\Sigma_{1} V_{1}}=\norm{\Sigma_{1}}$. In addition, the distribution of a Gaussian matrix is invariant under orthogonal transformations. Therefore, we can write:

\begin{flalign} \label{tr1}
	E(E(\norm{\Sigma_{3}\mathscr{G}_{3} \hat{\mathscr{G}_{2}}}^{2})|\mathscr{G}_{2})  & =  E(E(\Sigma_{jk}(\sigma_{jj} [\mathscr{G}_{3}]_{jk} [\hat{\mathscr{G}_{2}}]_{kk}))\nonumber \\
	& = E(\Sigma_{jk}(\sigma^{2}_{jj} [\hat{\mathscr{G}_{2}}]^{2}_{kk})) \nonumber \\
	& = E (\norm{\Sigma_{3}}_{F}^{2} \norm{\hat{\mathscr{G}_{2}}}^{2}) \nonumber \\
	& = \norm{\Sigma_{3}}_{F}^{2} E(\norm{\hat{\mathscr{G}_{2}}}^{2}) \nonumber \\
	& = \frac{r}{s-1} \norm{\Sigma_{3}}_{F}^{2} \nonumber\\
	& = \frac{r}{s-1}(\sum_{i=r+1}^{\min{(a,b)}}\sigma_{i}^2) \nonumber\\
	\end{flalign}
	
	Therefore, putting the expression of $E(\norm{\Sigma_{3}\mathscr{G}_{3} \hat{\mathscr{G}_{2}}}^{2})$ in the Eq. \ref{f4}, we can write:
	\begin{equation}
E(\norm{F-H(H^{*}F)})) \leq (1+ \frac{r}{s-1}) ( \sum_{i=r+1}^{\min{(a,b)}} \sigma_{i}^2)
\label{f5}
\end{equation}

The Eckart–Young theorem \cite{eckart1936approximation} states that $\sum_{i=r+1}^{\min{(a,b)}} \sigma_{i}^2$ is the smallest possible error that can be incurred when approximating $F$ by a matrix of rank $r$. Therefore, in Eq. \ref{f5} the
optimal error bound is missed by a factor of $(1+ \frac{r}{s-1})$. Next our objective is to determine the value of $s$ for a given $\epsilon$.
\end{proof}
\subsection{Computing the value of $s$ in algorithm \ref{a1}} 
\begin{lemma}
$HT$ be an approximation of $F$ calculated by algorithm \ref{a1} and $s$ is an oversampling parameter. Then $\norm{F-HT}< \epsilon$ if $s= \ceil*{\frac{r(\sum_{i=r+1}^{\min{(a,b)}} \sigma_{i}^2)}{\epsilon-(\sum_{i=r+1}^{\min{(a,b)}} \sigma_{i}^2)}+1}$.
\label{l3}
\end{lemma}

\begin{proof}

The approximation of $F$ derived by algorithm \ref{a1}, should haver approximation error less than $\epsilon$. Eq. \ref{f5} gives us expected approximation error bound. Therefore from Eq. \ref{f5}, we can write:

\begin{flalign} 
	(1+ \frac{r}{s-1}) (\sum_{i=r+1}^{\min{(a,b)}} \sigma_{i}^2)  & < \epsilon \nonumber \\
	(1+ \frac{r}{s-1})  & < \frac{\epsilon}{(\sum_{i=r+1}^{\min{(a,b)}} \sigma_{i}^2)} \nonumber \\
	 \frac{r}{s-1}  & < \frac{\epsilon}{(\sum_{i=r+1}^{\min{(a,b)}} \sigma_{i}^2)} -1 \nonumber \\
	 \frac{r}{s-1}  & < \frac{\epsilon-(\sum_{i=r+1}^{\min{(a,b)}} \sigma_{i}^2)}{(\sum_{i=r+1}^{\min{(a,b)}} \sigma_{i}^2)}  \nonumber \\
	 \frac{r(\sum_{i=r+1}^{\min{(a,b)}} \sigma_{i}^2)}{\epsilon-(\sum_{i=r+1}^{\min{(a,b)}} \sigma_{i}^2)}  & < s-1  \nonumber \\
	 \frac{r(\sum_{i=r+1}^{\min{(a,b)}} \sigma_{i}^2)}{\epsilon-(\sum_{i=r+1}^{\min{(a,b)}} \sigma_{i}^2)}+1 & < s  \nonumber \\
	\end{flalign}

Therefore, we can choose $s= \ceil*{\frac{r(\sum_{i=r+1}^{\min{(a,b)}} \sigma_{i}^2)}{\epsilon-(\sum_{i=r+1}^{\min{(a,b)}} \sigma_{i}^2)}+1}$, where $\ceil*{.}$ function gives the least integer greater than or equal to the given input.
\end{proof}

\textbf{The upshot:} Putting the value of $s$ in algorithm \ref{a1} we can compute the $H$ and $T= H^{*}F$ such that $HT$ is an approximation of $F$ with approximation error less than $\epsilon$. Thus we prove that for a given $F$ and an arbitrary state-of-the-art approximation of $F$ say $F^{*}$, we can always find an approximation in factor form of $(HT)$ which is better than $F^{*}$. Recall that we
are working with the Frobenius norm of matrices, which seems to be the most common way to measure the cost function in the low-rank matrix approximation problem \cite{halko2011finding}\cite{martinsson2016randomized}. However, the same method will surely allow one to prove the result for low-rank matrix approximation with respect to other norms like the spectral norm of matrices.
\end{proof}

%\section*{References}
\section{Declaration of competing interest}
None declared.

\bibliography{elsarticle-template.bib}

\begin{thebibliography}{11}
\expandafter\ifx\csname natexlab\endcsname\relax\def\natexlab#1{#1}\fi
\providecommand{\url}[1]{\texttt{#1}}
\providecommand{\href}[2]{#2}
\providecommand{\path}[1]{#1}
\providecommand{\DOIprefix}{doi:}
\providecommand{\ArXivprefix}{arXiv:}
\providecommand{\URLprefix}{URL: }
\providecommand{\Pubmedprefix}{pmid:}
\providecommand{\doi}[1]{\href{http://dx.doi.org/#1}{\path{#1}}}
\providecommand{\Pubmed}[1]{\href{pmid:#1}{\path{#1}}}
\providecommand{\bibinfo}[2]{#2}
\ifx\xfnm\relax \def\xfnm[#1]{\unskip,\space#1}\fi
%Type = Article
\bibitem[{Achlioptas and McSherry(2007)}]{achlioptas2007fast}
\bibinfo{author}{D.~Achlioptas}, \bibinfo{author}{F.~McSherry},
\newblock \bibinfo{title}{Fast computation of low-rank matrix approximations},
\newblock \bibinfo{journal}{Journal of the ACM (JACM)} \bibinfo{volume}{54}
  (\bibinfo{year}{2007}) \bibinfo{pages}{9--es}.
%Type = Article
\bibitem[{Gittens and Tropp(2009)}]{gittens2009error}
\bibinfo{author}{A.~Gittens}, \bibinfo{author}{J.~A. Tropp},
\newblock \bibinfo{title}{Error bounds for random matrix approximation
  schemes},
\newblock \bibinfo{journal}{arXiv preprint arXiv:0911.4108}
  (\bibinfo{year}{2009}).
%Type = Article
\bibitem[{Drineas et~al.(2006)Drineas, Kannan, and Mahoney}]{drineas2006fast}
\bibinfo{author}{P.~Drineas}, \bibinfo{author}{R.~Kannan},
  \bibinfo{author}{M.~W. Mahoney},
\newblock \bibinfo{title}{Fast monte carlo algorithms for matrices ii:
  Computing a low-rank approximation to a matrix},
\newblock \bibinfo{journal}{SIAM Journal on computing} \bibinfo{volume}{36}
  (\bibinfo{year}{2006}) \bibinfo{pages}{158--183}.
%Type = Article
\bibitem[{Frieze et~al.(2004)Frieze, Kannan, and Vempala}]{frieze2004fast}
\bibinfo{author}{A.~Frieze}, \bibinfo{author}{R.~Kannan},
  \bibinfo{author}{S.~Vempala},
\newblock \bibinfo{title}{Fast monte-carlo algorithms for finding low-rank
  approximations},
\newblock \bibinfo{journal}{Journal of the ACM (JACM)} \bibinfo{volume}{51}
  (\bibinfo{year}{2004}) \bibinfo{pages}{1025--1041}.
%Type = Article
\bibitem[{Papadimitriou et~al.(2000)Papadimitriou, Raghavan, Tamaki, and
  Vempala}]{papadimitriou2000latent}
\bibinfo{author}{C.~H. Papadimitriou}, \bibinfo{author}{P.~Raghavan},
  \bibinfo{author}{H.~Tamaki}, \bibinfo{author}{S.~Vempala},
\newblock \bibinfo{title}{Latent semantic indexing: A probabilistic analysis},
\newblock \bibinfo{journal}{Journal of Computer and System Sciences}
  \bibinfo{volume}{61} (\bibinfo{year}{2000}) \bibinfo{pages}{217--235}.
%Type = Inproceedings
\bibitem[{Clarkson and Woodruff(2009)}]{clarkson2009numerical}
\bibinfo{author}{K.~L. Clarkson}, \bibinfo{author}{D.~P. Woodruff},
\newblock \bibinfo{title}{Numerical linear algebra in the streaming model},
\newblock in: \bibinfo{booktitle}{Proceedings of the forty-first annual ACM
  symposium on Theory of computing}, \bibinfo{year}{2009}, pp.
  \bibinfo{pages}{205--214}.
%Type = Article
\bibitem[{Mahoney and Drineas(2009)}]{mahoney2009cur}
\bibinfo{author}{M.~W. Mahoney}, \bibinfo{author}{P.~Drineas},
\newblock \bibinfo{title}{Cur matrix decompositions for improved data
  analysis},
\newblock \bibinfo{journal}{Proceedings of the National Academy of Sciences}
  \bibinfo{volume}{106} (\bibinfo{year}{2009}) \bibinfo{pages}{697--702}.
%Type = Article
\bibitem[{Goreinov et~al.(1997)Goreinov, Tyrtyshnikov, and
  Zamarashkin}]{goreinov1997theory}
\bibinfo{author}{S.~A. Goreinov}, \bibinfo{author}{E.~E. Tyrtyshnikov},
  \bibinfo{author}{N.~L. Zamarashkin},
\newblock \bibinfo{title}{A theory of pseudoskeleton approximations},
\newblock \bibinfo{journal}{Linear algebra and its applications}
  \bibinfo{volume}{261} (\bibinfo{year}{1997}) \bibinfo{pages}{1--21}.
%Type = Article
\bibitem[{Halko et~al.(2011)Halko, Martinsson, and Tropp}]{halko2011finding}
\bibinfo{author}{N.~Halko}, \bibinfo{author}{P.-G. Martinsson},
  \bibinfo{author}{J.~A. Tropp},
\newblock \bibinfo{title}{Finding structure with randomness: Probabilistic
  algorithms for constructing approximate matrix decompositions},
\newblock \bibinfo{journal}{SIAM review} \bibinfo{volume}{53}
  (\bibinfo{year}{2011}) \bibinfo{pages}{217--288}.
%Type = Article
\bibitem[{Martinsson and Voronin(2016)}]{martinsson2016randomized}
\bibinfo{author}{P.-G. Martinsson}, \bibinfo{author}{S.~Voronin},
\newblock \bibinfo{title}{A randomized blocked algorithm for efficiently
  computing rank-revealing factorizations of matrices},
\newblock \bibinfo{journal}{SIAM Journal on Scientific Computing}
  \bibinfo{volume}{38} (\bibinfo{year}{2016}) \bibinfo{pages}{S485--S507}.
%Type = Article
\bibitem[{Eckart and Young(1936)}]{eckart1936approximation}
\bibinfo{author}{C.~Eckart}, \bibinfo{author}{G.~Young},
\newblock \bibinfo{title}{The approximation of one matrix by another of lower
  rank},
\newblock \bibinfo{journal}{Psychometrika} \bibinfo{volume}{1}
  (\bibinfo{year}{1936}) \bibinfo{pages}{211--218}.

\end{thebibliography}

\end{document}